\documentclass[11pt]{amsart}
\usepackage{fullpage}
\usepackage[utf8]{inputenc}
\usepackage[T1]{fontenc}
\usepackage{graphicx}
\usepackage{amsmath}
\usepackage{amsfonts}
\usepackage{graphicx}

\newtheorem{theo}{Theorem}[section]
\newtheorem{cor}[theo]{Corollary}
\newtheorem{rem}[theo]{Remark}

\newtheorem{propo}[theo]{Proposition}
\newtheorem{lemme}[theo]{Lemma}
\newtheorem{defi}[theo]{Definition}
\newtheorem{ex}[theo]{Example}

\newtheorem{hyp}[theo]{Assumptions}
\newcommand{\E}{\mathbb{E}}
\newcommand{\R}{\mathbb{R}}
\newcommand{\PP}{\mathbb{P}}

\newcommand{\N}{\mathbb{N}}

\numberwithin{equation}{section}

\title{Analysis of a HMM time-discretization scheme for a system of Stochastic PDE's}
\author{Charles-Edouard BREHIER}
\date{}
\address{ENS Cachan Bretagne - IRMAR, Universit\'e Rennes 1\\
Avenue Robert Schumann\\ F-35170 Bruz\\ France}

\email{charles-edouard.brehier@bretagne.ens-cachan.fr}

\keywords{Stochastic Partial Differential Equations, Heterogeneous Multiscale Method, Time Scale Separation, Euler Scheme, Strong and Weak Approximation}
\subjclass{60H15,60H35,70K70}

\begin{document}
\maketitle
\begin{abstract}
We consider the discretization in time of a system of parabolic stochastic partial differential equations with slow and fast components; the fast equation is driven by an additive space-time white noise. The numerical method is inspired by the Averaging Principle satisfied by this system, and fits to the framework of Heterogeneous Multiscale Methods.

The slow and the fast components are approximated with two coupled numerical semi-implicit Euler schemes depending on two different timestep sizes.

We derive bounds of the approximation error on the slow component in the strong sense - approximation of trajectories - and in the weak sense - approximation of the laws. The estimates generalize the results of \cite{E-L-V} in the case of infinite dimensional processes.

\end{abstract}




\section{Introduction}

In this paper, we are interested in the numerical approximation of a randomly-perturbed system of reaction-diffusion equations that can be written
\begin{equation}\label{truePDE}
\begin{gathered}
\frac{\partial x^\epsilon(t,\xi)}{\partial t}=\frac{\partial^2 x^\epsilon(t,\xi)}{\partial \xi^2}+f(\xi,x^\epsilon(t,\xi),y^\epsilon(t,\xi)),\\
\frac{\partial y^\epsilon(t,\xi)}{\partial t}=\frac{1}{\epsilon}\frac{\partial^2 y^\epsilon(t,\xi)}{\partial \xi^2}+\frac{1}{\epsilon}g(\xi,x^\epsilon(t,\xi),y^\epsilon(t,\xi))+\frac{1}{\sqrt{\epsilon}}\frac{\partial \omega(t,\xi)}{\partial t},
\end{gathered}
\end{equation}
for $t\geq 0, \xi\in(0,1)$, with initial conditions $x^\epsilon(0,\xi)=x(\xi)$ and $y^\epsilon(0,\xi)=\xi$, and homogeneous Dirichlet boundary conditions $x^\epsilon(t,0)=x^\epsilon(t,1)=0,y^\epsilon(t,0)=y^\epsilon(t,1)=0$.
The stochastic perturbation $\frac{\partial \omega(t,\xi)}{\partial t}$ is a space-time white noise and $\epsilon>0$ is a small parameter.


In the recent article \cite{CEB1}, we have proved that an Averaging Principle holds for such a system, and we have exhibited an order of convergence - with respect to $\epsilon$ - in a strong and in a weak sense: the slow component $x^\epsilon$ is approximated thanks to the solution of an averaged equation. In this article, we analyse a numerical method of time discretization which reproduces this averaging effect at the discrete time level. More precisely, our aim is to build a numerical approximation of the slow component $x^\epsilon$, taking care of the stiffness induced by the time scale separation. The Heterogeneous Multiscale Method - HMM - procedure can be used, as it is done in \cite{E-L-V} for SDEs of the same kind. First we recall the general principle of such a method, which has been developped in various contexts, both deterministic and stochastic - see the review article \cite{Review} and the references therein, as well as \cite{E1}, \cite{E-Eng}, \cite{Van1}.

In system \eqref{truePDE}, the two components evolve at different time scales; $x^\epsilon$ is the slow component of the problem, while the fast component $y^\epsilon$ has fast variations in time. We are indeed interested in evaluating the slow component, which can be thought as the mathematical model for a phenomenon appearing at the natural time scale of the experiment, whereas the fast component can often be interpreted as data for the slow component, taking care of effects at a faster time scale. Instead of using a direct numerical method, which might require a very small time step size because of the fast component, we use a different solver for each time scale: a macrosolver and a microsolver. The macrosolver leads to the approximation of the slow component; it takes into account data from the evolution at the fast time scale. The microsolver is then a procedure for estimating the unknown necessary data, using the evolution at the microtime scale, which also depends on the evolution at the macrotime scale. We emphasize a difference with the framework of \cite{E-L-V}: instead of analysing various numerical schemes, the infinite dimensional setting implies that we only focus on a semi-implicit Euler scheme.


In Section \ref{secttheo}, we state the two main Theorems of this article: we show a strong convergence result - Theorem \ref{cvforte} - as well as a weak convergence result - Theorem \ref{cvfaible} - which are similar to the available results for SDEs. Compared to \cite{E-L-V}, we propose modified and simplified proofs leading to apparently weaker error estimates; we made this choice for various reasons. First, even if apparently we get weaker estimates, under an appropriate choice of the parameters the cost of the method remains of the same order. Second, the generalization of the final dimensional results would not yield the same bounds, due to the regularity assumptions we make on the nonlinear coefficients of the equations. Finally, we can extend the weak convergence result to the situation where the fast equation only satisfies a weak dissipativity assumption.



In the case of a linear fast equation - when $g$ is equal to $0$ - it is well-known that the second equation in \eqref{truePDE} is dissipative. In the general case, we make assumptions on $g$ so that this property is preserved for $y^\epsilon$ - see Assumptions \ref{strictdiss} and \ref{weakdiss} below. The fast equation with frozen slow component - defined by \eqref{eqfig} in the abstract framework - then admits a unique invariant probability measure, which is ergodic and strongly mixing - with exponential convergence to equilibrium. Under the strict dissipativity condition \ref{strictdiss}, we can prove that the averaging principle holds in the strong and in the weak sense; moreover the ``fast'' numerical scheme has the same asymptotic behaviour as the continuous time equation. If we only assume weak dissipativity of Assumption \ref{weakdiss}, the averaging principle only holds in the weak sense, and we can not prove uniqueness of the invariant law of the fast numerical scheme. Nevertheless, in the general setting \cite{CEB2} gives an approximation result of the invariant law of the continuous time equation with the numerical method which is used to prove Theorem \ref{cvfaible}; the order of convergence is $1/2$ with respect to the timestep size - the precise result is recalled in Theorem \ref{weakdistance}.


The paper is organized as follows. In Section \ref{sectdescri}, we give the definition of the numerical scheme. We then state the main assumptions made on the system of equations. In Section \ref{secttheo} we state the two main Theorems proved in this article, while in Section \ref{commsect} we compare the efficiency of the HMM scheme with a direct one in order to justify the use of a new method. Before proving the theorems, we give some useful results on the numerical schemes. 
Finally the last two sections contain the proof of the strong and weak convergence theorems.


\section{Description of the numerical scheme}\label{sectdescri}

Instead of working directly with system \eqref{truePDE}, we work with abstract stochastic evolution equations in Hilbert spaces $H$:
\begin{equation}\label{eqgen}
\begin{gathered}
dX^{\epsilon}(t)=(AX^{\epsilon}(t)+F(X^{\epsilon}(t),Y^{\epsilon}(t)))dt\\
dY^{\epsilon}(t)=\frac{1}{\epsilon}(BY^{\epsilon}(t)+G(X^{\epsilon}(t),Y^{\epsilon}(t)))dt+\frac{1}{\sqrt{\epsilon}}dW(t),
\end{gathered}
\end{equation}
with initial conditions $X^{\epsilon}(0)=x\in H$, $Y^{\epsilon}(0)=y\in H$.

To get system \eqref{eqgen} from \eqref{truePDE}, we take $H=L^2(0,1)$; the linear operators $A$ and $B$ are Laplace operators with homogeneous Dirichlet boundary conditions - see Example \ref{exampleAB} - and the nonlinearities $F$ and $G$ are Nemytskii operators - see Example \ref{exFG}. The process $(W(t))_{t\geq 0}$ is a cylindrical Wiener process on $H$ - see Section \ref{SectWiener}. For precise assumptions on the coefficients, we refer to Section \ref{Sectassumptions}.

We recall the idea of the Averaging Principle - proved in the previous article \cite{CEB1}: when $\epsilon$ goes to $0$, $X^\epsilon$ can be approximated by $\overline{X}$ defined by the averaged equation
\begin{equation}\label{eqmoy}
\begin{gathered}
\frac{d\overline{X}(t)}{dt}=A\overline{X}(t)+\overline{F}(\overline{X}(t))\\
\overline{X}(t=0)=x\in H;
\end{gathered}
\end{equation}
the error is controlled in a strong sense by $C\epsilon^{1/2-r}$ and in a weak sense by $C\epsilon^{1-r}$, where $r>0$ can be chosen as small as necessary, and where $C$ is a constant.

The averaged coefficient $\overline{F}$ - see \eqref{deffbar} - satisfies
$$\overline{F}(X)=\int_{H}F(X,y)\mu^X(dy)=\lim_{s \to +\infty}\E[F(X,Y_X(s,y))],$$
where $\mu^X$ is the unique invariant probability measure of the fast process $Y_X$ with frozen slow component - more details are given in Section \ref{known}. 

To apply the HMM strategy, we need to define a macrosolver and a microsolver. We denote by $\Delta t$ the macrotime step size, and by $\delta t$ the microtime step size. Let also $T>0$ be a given final time.

The construction of the macrosolver is deeply based on the averaging principle: for $n\Delta t\leq T$ $X^\epsilon(n\Delta t)$ can be approximated by $\overline{X}(n\Delta t)$. If the averaged coefficient $\overline{F}$ was known, one could build an approximation with a deterministic numerical scheme on the averaged equation; nevertheless in general it is not the case, and the idea is to calculate an approximation of this coefficient on-the-fly, by using the microsolver.

Therefore the macrosolver is defined in the following way: for any $0\leq n\leq \lfloor \frac{T}{\Delta t} \rfloor:=n_0$,
$$X_{n+1}=X_{n}+\Delta t AX_{n+1}+\Delta t \tilde{F}_n,$$
with the initial condition $X_0=x$. $\tilde{F}_n$ has to be defined; before that, we notice that the above definition leads to a semi-implicit Euler scheme - we use implicitness on the linear part, but the nonlinear part is explicit. If we define a bounded linear operator on $H$ by $S_{\Delta t}=(I-\Delta t A)^{-1}$, we rather use the following explicit formula
\begin{equation}\label{schémalent}
X_{n+1}=S_{\Delta t}X_{n}+\Delta t S_{\Delta t}\tilde{F}_n.
\end{equation}

We want $\tilde{F}_n$ to be an approximation of $\overline{F}(X_n)$. The role of the microsolver is to give an approximation of $Y_{X_n}$ - the fast process with frozen slow component $X_n$, when $n$ is fixed; moreover we compute a finite number $M$ of independent replicas of the process, in order to approximate theoretical expectations by discrete averages over different realizations of the random variables, in a Monte-Carlo spirit. Therefore the microsolver is defined in the following way: we fix a realization index $j\in \left\{ 1,\ldots,M\right\}$, and a macrotime step $n$; then for any $m\geq 0$
$$Y_{n,m+1,j}=Y_{n,m,j}+\frac{\delta t}{\epsilon}BY_{n,m+1,j}+\frac{\delta t}{\epsilon}G(X_n,Y_{n,m,j})+\sqrt{\frac{\delta t}{\epsilon}}\zeta_{n,m+1,j}.$$
As above we can give an explicit formula
\begin{equation}\label{schémarapide}
Y_{n,m+1,j}=R_{\frac{\delta t}{\epsilon}}Y_{n,m,j}+\frac{\delta t}{\epsilon}R_{\frac{\delta t}{\epsilon}}G(X_n,Y_{n,m,j})+\sqrt{\frac{\delta t}{\epsilon}}R_{\frac{\delta t}{\epsilon}}\zeta_{n,m+1,j},
\end{equation}
where for any $\tau>0$ $R_{\tau}=(I-\tau B)^{-1}$.

The noises $\zeta_{n,m,j}$ are defined by
$$\zeta_{n,m+1,j}=\frac{W_{(m+1)\delta t}^{(n,j)}-W_{m\delta t}^{(n,j)}}{\sqrt{\delta t}},$$
where $(W^{(n,j)})_{1\leq j\leq M,0\leq n\leq n_0}$ are independent cylindrical Wiener processes on $H$. It is essential to use independent noises at each macrotime step. It is important to remark that this equation is well-posed in the Hilbert space $H$ - see the conditions required when a cylindrical Wiener process is used in Section \ref{SectWiener} - since $R_{\tau}$ is a Hilbert-Schmidt operator from $H$ to $H$, under assumptions given in Section \ref{Sectassumptions}.

The missing definition can now be written: $\tilde{F}_n$ is given by
\begin{equation}\label{ftilde}
\tilde{F}_n=\frac{1}{MN}\sum_{j=1}^{M}\sum_{m=n_T}^{n_T+N-1}F(X_n,Y_{n,m,j}).
\end{equation}
$n_T$ is the number of microtime steps that are not used in the evaluation of the average in $\tilde{F}_n$, while $N$ is the number of microtime steps that are then used for this evaluation. Each macrotime step then requires the computation of $m_0=n_T+N-1$ values of the microsolver.

At each macrotime step $Y_{n,m,j}$ must be initialized at time $m=0$. In our proofs, this is not as important as in \cite{E-L-V}, but for definiteness we use the same method: we initialize with the last value computed during the previous macrotime step: $Y_{n+1,0,j}=Y_{n,m_0,j}$, while $Y_{0,0,j}=y$.



The aim of the analysis for HMM schemes is to prove that under an appropriate choice of the parameters $n_T,N,M$ of the scheme, we can bound the error by expressions of the following kind, where $r>0$ is chosen as small as necessary: for $n=n_0=\lfloor \frac{T}{\Delta t} \rfloor$, we have the strong error estimate
$$\E|X^\epsilon(n\Delta t)-X_n|\leq C\left(\epsilon^{1/2-r}+\Delta t^{1-r}+\left(\frac{\delta t}{\epsilon}\right)^{1/2-r}\right),$$
and if $\Phi$ is a test function of class $\mathcal{C}^2$, bounded and with bounded derivatives, we have the weak error estimate
$$|\E\Phi(X^\epsilon(n\Delta t))-\E\Phi(X_n)|\leq C\left(\epsilon^{1-r}+\Delta t^{1-r}+\left(\frac{\delta t}{\epsilon}\right)^{1/2-r}\right).$$
The origin of the three error terms appears clearly in the proofs - see Sections \ref{sectfort} and \ref{sectfaible}: the first one is the averaging error, the second one is the error in a deterministic scheme with the macrotime step, and the third one is the weak error in a scheme for stochastic equations with the microtime step. We recall that in the SPDE case the strong order of the semi-implicit Euler scheme for the microsolver used here is $1/4$, while the weak order is $1/2$, while in the SDE situation the respective orders are $1/2$ and $1$. The macrosolver is deterministic, so that the order is $1$.
Precise results for any choice of $n_T,N,M$ are given in Theorems \ref{cvforte} and \ref{cvfaible} below, while the choice of these parameters is explained in Section \ref{commsect}.

\section{Assumptions}\label{Sectassumptions}

As mentioned above, System \eqref{eqgen} satisfies an averaging principle, and strong and weak order of convergence can be given. The HMM method relies on that idea. The natural assumptions are basically the same as the hypothesis needed to prove those results, but must be strenghtened sometimes. The typical kind of coefficients is specified in Examples \ref{exampleAB} and \ref{exFG}.

\subsection{The cylindrical Wiener process and stochastic integration in H}\label{SectWiener}

Here, we recall the definition of the cylindrical Wiener process and of stochastic integral on a separable Hilbert space $H$ - its norm is denoted by $|.|_{H}$ or just $|.|$. For more details, see \cite{DaP-Z1}.

We first fix a filtered probability space $(\Omega,\mathcal{F},(\mathcal{F}_t)_{t\geq 0},\PP)$. A cylindrical Wiener process on $H$ is defined with two elements:
\begin{itemize}
\item a complete orthonormal system of $H$, denoted by  $(q_i)_{i\in I}$, where $I$ is a subset of $\N$;
\item a family $(\beta_i)_{i\in I}$ of independent real Wiener processes with respect to the filtration $((\mathcal{F}_t)_{t\geq 0})$:
\end{itemize}

\begin{equation}\label{defWiener}
W(t)=\sum_{i\in I}\beta_i(t)q_i.
\end{equation}

When $I$ is a finite set, we recover the usual definition of Wiener processes in the finite dimensional space $\R^{|I|}$. However the subject here is the study of some Stochastic Partial Differential Equations, so that in the sequel the underlying Hilbert space $H$ is infinite dimensional; for instance when $H=L^{2}(0,1)$, an example of complete orthonormal system is $(q_k)=(\sin(k.))_{k\geq 1}$ - see Example \ref{exampleAB}.

A fundamental remark is that the series in \eqref{defWiener} does not converge in $H$; but if a linear operator $\Psi:H\rightarrow K$ is Hilbert-Schmidt, then $\Psi W(t)$ converges in $L^2(\Omega,H)$ for any $t\geq 0$.

We recall that a linear operator $\Psi:H\rightarrow K$ is said to be Hilbert-Schmidt when
$$|\Psi|_{\mathcal{L}_{2}(H,K)}^{2}:=\sum_{k=0}^{+\infty}|\Psi(q_k)|_{K}^{2}<+\infty,$$
where the definition is independent of the choice of the orthonormal basis $(q_k)$ of $H$.
The space of Hilbert-Schmidt operators from $H$ to $K$ is denoted $\mathcal{L}_{2}(H,K)$; endowed with the norm $|.|_{\mathcal{L}_{2}(H,K)}$ it is an Hilbert space.

The stochastic integral $\int_{0}^{t}\Psi(s)dW(s)$ is defined in $K$ for predictible processes $\Psi$ with values in $\mathcal{L}_2(H,K)$ such that $\int_{0}^{t}|\Psi(s)|_{\mathcal{L}_2(H,K)}^{2}ds<+\infty$ a.s; moreover when $\Psi\in L^2(\Omega\times[0,t];\mathcal{L}_2(H,K))$, the following two properties hold:
\begin{gather*}
\E|\int_{0}^{t}\Psi(s)dW(s)|_{K}^{2}=\E\int_{0}^{t}|\Psi(s)|_{\mathcal{L}_2(H,K)}^{2}ds, \text{ (It\^o isometry),}\\
\E\int_{0}^{t}\Psi(s)dW(s)=0.
\end{gather*}
A generalization of It\^o formula also holds - see \cite{DaP-Z1}.

For instance, if $v=\sum_{k\in\N}v_kq_k\in H$, we can define $$<W(t),v>=\int_{0}^{t}<v,dW(s)>=\sum_{k\in\N}\beta_k(t)v_k;$$
we then have the following space-time white noise property
$$\E<W(t),v_1><W(s),v_2>=t\wedge s<v_1,v_2>.$$

Therefore to be able to integrate a process with respect to $W$ requires some strong properties on the integrand; in our SPDE setting, the Hilbert-Schmidt properties follow from the assumptions made on the linear coefficients of the equations.

\subsection{Assumptions on the linear operators}
First, we have to specify some properties of the linear operators $A$ and $B$ coming into the definition of system \eqref{eqgen}; we assume that the linear parts are of parabolic type, with space variable $\xi\in(0,1)$.

We assume that $A$ and $B$ are unbounded linear operators, with domains $D(A)$ and $D(B)$, which satisfy the following assumptions:
\begin{hyp}\label{hypAB}
\begin{enumerate}
\item There exist complete orthonormal systems of $H$, $(e_k)_{k\in \N}$ and $(f_k)_{k\in \N}$, and $(\lambda_k)_{k\in \N}$ and $(\mu_k)_{k\in \N}$ non-decreasing sequences of real positive numbers such that:
\begin{gather*}
Ae_{k}=-\lambda_{k}e_{k}\text{ for all } k\in\N\\
Bf_{k}=-\mu_{k}f_{k}\text{ for all } k\in\N.
\end{gather*}
We use the notations $\lambda:=\lambda_0>0$ and $\mu:=\mu_0>0$ for the smallest eigenvalues of $A$ and $B$.
\item For every $k\in \N$, $f_k$ is Lipschitz continuous and bounded on $[0,1]$, with a uniform control with respect to $k$: there exists $C>0$ such that for any $\xi_1,\xi_2\in [0,1]$
$$|f_k(\xi_1)|\leq C \quad \text{ and } \quad |f_k(\xi_1)-f_k(\xi_2)|\leq C\sqrt{\mu_k}|\xi_1-\xi_2|.$$
\item The sequences $(\lambda_k)$ and $(\mu_k)$ go to $+\infty$; moreover we have some control of the behaviour of $(\mu_k)$ given by:
$$
\sum_{k=0}^{+\infty}\frac{1}{\mu_{k}^{\alpha}}<+\infty\Leftrightarrow \alpha>1/2.$$
\end{enumerate}
\end{hyp}

\begin{ex}\label{exampleAB}
$A=B=\frac{d^2}{dx^2}$, with domain $H^2(0,1)\cap H_{0}^{1}(0,1)\subset L^2(0,1)$ - homogeneous Dirichlet boundary conditions: in that case $\lambda_k=\mu_k=\pi^2 k^2$, and $e_k(\xi)=f_k(\xi)=\sqrt{2}\sin(k\pi\xi)$ - see \cite{Brezis}.
\end{ex}

In the abstract setting, powers of $-A$ and $-B$, with their domains can be easily defined:

\begin{defi}\label{defpowers}
For $a,b\in[0,1]$, we define the operators $(-A)^a$ and $(-B)^b$ by
\begin{gather*}
(-A)^a x=\sum_{k=0}^{\infty}\lambda_{k}^{a}x_ke_k\in H,\\
(-B)^b y=\sum_{k=0}^{\infty}\mu_{k}^{b}y_kf_k\in H,
\end{gather*}

with domains
\begin{gather*}
D(-A)^a=\left\{x=\sum_{k=0}^{+\infty}x_ke_{k}\in H; |x|_{(-A)^a}^{2}:=\sum_{k=0}^{+\infty}(\lambda_k)^{2a}|x_k|^2<+\infty\right\};\\
D(-B)^b=\left\{y=\sum_{k=0}^{+\infty}y_kf_{k}\in H, |y|_{(-B)^b}^{2}:=\sum_{k=0}^{+\infty}(\mu_k)^{2b}|y_k|^2<+\infty\right\}.
\end{gather*}
\end{defi}

On $D(-A)^a$, the norm $|\hspace{3pt}.\hspace{3pt}|_{(-A)^a}$ and the Sobolev norm of $H^{2a}$ are equivalent: when $x$ belongs to a space $D(-A)^a$, the exponent $a$ represents some regularity of the function $x$.

When $a\geq 0$, we can also define a bounded linear operator $(-A)^{-a}$ in $H$ with
$$(-A)^{-a}x=\sum_{k=0}^{+\infty}\lambda_{k}^{-a}x_k\in H,$$
where $x=\sum_{k=0}^{+\infty}x_ke_{k}\in H$.

Under the previous assumptions on the linear coefficients, it is easy to show that the following stochastic integral is well-defined in $H$, for any $t\geq 0$:
\begin{equation}\label{stoconvWB}
W^B(t)=\int_{0}^{t}e^{(t-s)B}dW(s).
\end{equation}

It is called a stochastic convolution, and it is the unique mild solution of
$$dZ(t)=BZ(t)dt+dW(t), \quad Z(0)=0.$$

Under the second condition of Assumption \ref{hypAB}, there exists $\delta>0$ such that for any $t>0$ we have $\int_{0}^{t}\frac{1}{s^\delta}|e^{sB}|_{\mathcal{L}_{2}(H)}^{2}ds<+\infty$; it can then be proved that $W^B$ has continuous trajectories - via the \textit{factorization method}, see \cite{DaP-Z1} - and that for any $1\leq p<+\infty$, any $0\leq \gamma<1/4$, there exists a constant $C(p,\gamma)<+\infty$ such that for any $t\geq 0$
\begin{equation}\label{estimWB}
\E|W^B(t)|_{(-B)^\gamma}^{p}\leq C(p,\gamma).
\end{equation}


\subsection{Assumptions on the nonlinear coefficients}
We now give the Assumptions on the nonlinear coefficients $F,G:H\times H\rightarrow H$. First, we need some regularity properties:
\begin{hyp}\label{hypF}
We assume that there exists $0\leq \eta<\frac{1}{2}$ and a constant $C$ such that the following directional derivatives are well-defined and controlled:
\begin{itemize}
\item For any $x,y\in H$ and $h\in H$, $|D_xF(x,y).h|\leq C|h|_{H}$ and $|D_yF(x,y).h|\leq C|k|_{H}$.
\item For any $x,y\in H$, $h\in H$, $k\in D(-A)^{\eta}$, $|D_{xx}^{2}F(x,y).(h,k)|\leq C|h|_{H}|k|_{(-A)^{\eta}}$.
\item For any $x,y\in H$, $h\in H$, $k\in D(-B)^{\eta}$, $|D_{yy}^{2}F(x,y).(h,k)|\leq C|h|_{H}|k|_{(-B)^{\eta}}$.
\item For any $x,y\in H$, $h\in H$, $k\in D(-B)^{\eta}$, $|D_{xy}^{2}F(x,y).(h,k)|\leq C|h|_{H}|k|_{(-B)^{\eta}}$.
\item For any $x,y\in H$, $h\in D(-A)^{\eta}$, $k\in H$, $|D_{xy}^{2}F(x,y).(h,k)|\leq C|h|_{(-A)^{\eta}}|k|_{H}$.
\end{itemize}
We moreover assume that $F$ is bounded.
\end{hyp}

We also need:
\begin{hyp}\label{hypF2}
For $\eta$ defined in the previous Assumption \ref{hypF}, we have for any $x,y\in H$ and $h,k\in H$
\begin{gather*}
|(-A)^{-\eta}D_{xx}^{2}F(x,y).(h,k)|\leq C|h|_H|k|_H\\
|(-A)^{-\eta}D_{xy}^{2}F(x,y).(h,k)|\leq C|h|_H|k|_H\\
|(-A)^{-\eta}D_{yy}^{2}F(x,y).(h,k)|\leq C|h|_H|k|_H.
\end{gather*}
\end{hyp}

We assume that the fast equation is a gradient system: for any $x$ the nonlinear coefficient $G(x,.)$ is the derivative of some potential $U$. We also assume regularity assumptions as for $F$.
\begin{hyp}\label{hypG}
The function $G$ is defined through $G(x,y)=\nabla_yU(x,y)$, for some potential $U:H\times H\rightarrow \R$. Moreover we assume that $G$ is bounded, and that the regularity assumptions given in the Assumption \ref{hypF} are also satisfied for $G$.
\end{hyp}

For $G$, we need a stronger hypothesis than for $F$ - in order to get Proposition \ref{FbarLip} . Assumption \ref{hypF2} becomes:
\begin{hyp}\label{hypG2}
We have for any $x,y\in H$, $h,k\in H$, $z\in L^\infty(0,1)$
$$|<D_{xx}G(x,y).(h,k),z>_{H}|\leq C|h|_H|k|_H|z|_{L^\infty(0,1)}.$$
\end{hyp}

Finally, we need to assume some dissipativity of the fast equation. Assumption \ref{strictdiss} is necessary to obtain strong convergence in the Averaging Principle, while Assumption \ref{weakdiss} is weaker and can lead to the weak convergence.

\begin{hyp}[Strict dissipativity]\label{strictdiss}
Let $L_g$ denote the Lipschitz constant of $G$ with respect to its second variable; then
\begin{equation}\label{hypLg}\tag{SD}
L_{g}<\mu,
\end{equation}
where $\mu$ is defined in Assumption \ref{hypAB}.
\end{hyp}

\begin{hyp}[Weak Dissipativity]\label{weakdiss}
There exist $c>0$ and $C>0$ such that for any $y\in D(B)$
\begin{equation}\label{hypdiss}\tag{WD}
<By+G(y),y>\leq -c|y|^2+C.
\end{equation}
\end{hyp}

The second Assumption is satisfied as soon as $G$ is bounded, while the first one requires some knowledge of the Lipschitz constant of $G$.

\begin{ex}\label{exFG}
We give some fundamental examples of nonlinearities for which the previous assumptions are satisfied:
\begin{itemize}
\item Functions $F,G:H\times H\rightarrow H$ of class $\mathcal{C}^2$, bounded and with bounded derivatives, such that $G(x,y)=\nabla_yU(x,y)$ and satisfying \eqref{hypLg} fit in the framework, with the choice $\eta=0$.
\item Functions $F$ and $G$ can be \textbf{Nemytskii} operators: let $f:(0,1)\times \R^2\rightarrow \R$ be a measurable function such that for almost every $\xi\in(0,1)$ $f(\xi,.)$ is twice continuously differentiable, bounded and with uniformly bounded derivatives.
Then $F$ is defined for every $x,y\in H=L^2(0,1)$ by
$$F(x,y)(\xi)=f(\xi,x(\xi),y(\xi)).$$
For $G$, we assume that there exists a function $g$ with the same properties as $f$ above, such that $G(x,y)(\xi)=g(\xi,x(\xi),y(\xi))$. The strict dissipativity Assumption \ref{hypLg} is then satisfied when $$\sup_{\xi\in(a,b),x\in \R,y\in \R}|\frac{\partial g}{\partial y}(\xi,x,y)|<\mu.$$
The conditions in Assumption \ref{hypF} are then satisfied for $F$ and $G$ as soon as there exists $\eta<1/2$ such that $D(-A)^\eta$ and $D(-B)^\eta$ are continuously embedded into $L^\infty(0,1)$ - it is the case for $A$ and $B$ given in Example \ref{exampleAB}, with $\eta>1/4$. Assumptions \ref{hypF2} and \ref{hypG2} are also satisfied.
\end{itemize}
\end{ex}

We then deduce that the system \eqref{eqgen} is well-posed for any $\epsilon>0$, on any finite time interval $[0,T]$. Under Assumptions \ref{hypAB}, \ref{hypF}, \ref{hypG}, the nonlinearities $F$ and $G$ are Lipschitz continuous, and the following Proposition is classical - see \cite{DaP-Z1}:
\begin{propo}
For every $\epsilon>0$, $T>0$, $x\in H$, $y\in H$, system \eqref{eqgen} admits a unique mild solution $(X^{\epsilon},Y^{\epsilon})\in (\text{L}^2(\Omega,\mathcal{C}([0,T],H)))^2$: for any $t\in[0,T]$,
\begin{equation}\label{eqmild}
\begin{gathered}
X^{\epsilon}(t)=e^{tA}x+\int_{0}^{t}e^{(t-s)A}F(X^{\epsilon}(s),Y^{\epsilon}(s))ds\\
Y^{\epsilon}(t)=e^{\frac{t}{\epsilon}B}y+\frac{1}{\epsilon}\int_{0}^{t}e^{\frac{(t-s)}{\epsilon}B}G(X^{\epsilon}(s),Y^{\epsilon}(s))ds
+\frac{1}{\sqrt{\epsilon}}\int_{0}^{t}e^{\frac{(t-s)}{\epsilon}B}dW(s).
\end{gathered}
\end{equation}
\end{propo}

\section{Convergence results}

\subsection{Statement of the Theorems}\label{secttheo}




We can now state the main results: the numerical process $(X_n)$ defined by \eqref{schémalent} approximates the slow component $X^\epsilon(n\Delta t)$ of system \eqref{eqgen}, with strong and weak error estimates given in the Theorems.

If the Strict Dissipativity Assumption \ref{strictdiss} is satisfied, we can prove the following:
\begin{theo}[Strong convergence]\label{cvforte}
Assuming \eqref{hypLg}, for any $0<r<1/2$, $0<\kappa<1/2$, $T>0$, $\epsilon_0>0$, $\Delta_0>0$, $\tau_0>0$, there exists $C>0$ such that for any $0<\epsilon\leq \epsilon_0$, $0<\Delta t\leq \Delta_0$, $\delta t>0$ such that $\tau=\frac{\delta t}{\epsilon}\leq \tau_0$ and $1<n\leq n_0=\lfloor\frac{T}{\Delta t}\rfloor$
\begin{align*}
\E|X^\epsilon(n\Delta t)-X_n|&\leq C\biggl(\epsilon^{1/2-r}+\frac{1}{n}+\Delta t^{1-r}\biggr)\\
&+C\Biggl((\frac{\delta t}{\epsilon})^{1/2-\kappa}+\frac{1}{\sqrt{N\frac{\delta t}{\epsilon}+1}}e^{-cn_T\frac{\delta t}{\epsilon}}\Biggr)\\
&+C\frac{\sqrt{\Delta t}}{\sqrt{M(N\frac{\delta t}{\epsilon}+1)}}.
\end{align*}


\end{theo}

Under the more general Weak Dissipativity Assumption \ref{weakdiss}, we can prove the following:
\begin{theo}[Weak convergence]\label{cvfaible}
Assume that $x\in D((-A)^\theta)$ and $y\in D((-B)^\beta)$, with $\theta,\beta\in]0,1]$.
Let $\Phi:H\rightarrow H$ bounded, of class $\mathcal{C}^2$, with bounded first and second order derivatives.
Then with the weak dissipativity assumption \eqref{hypdiss}, for any $0<r<1$, $\kappa<1/2$, $T>0$, $\epsilon_0>0$, $\Delta_0>0$, $\tau_0>0$, there exists $C>0$ such that for any $0<\epsilon\leq \epsilon_0$, $0<\Delta t\leq \Delta_0$, $\delta t>0$ such that $\tau=\frac{\delta t}{\epsilon}\leq \tau_0$ and $1<n\leq n_0=\lfloor\frac{T}{\Delta t}\rfloor$
\begin{align*}
|\E\Phi(X^\epsilon(n\Delta t))-\E\Phi(X_n)|&\leq C\biggl(\epsilon^{1-r}+\frac{1}{n}+\Delta t^{1-r}\biggr)\\
&+C\Biggl((\frac{\delta t}{\epsilon})^{1/2-\kappa}(1+\frac{1}{((n_T-1)\frac{\delta t}{\epsilon})^{1/2-\kappa}})+\frac{1}{N\frac{\delta t}{\epsilon}+1}e^{-cn_T\frac{\delta t}{\epsilon}}\Biggr).
\end{align*}
\end{theo}

Moreover, if we assume \eqref{hypLg} in Theorem \ref{cvfaible}, the factor $(1+\frac{1}{((n_T-1)\frac{\delta t}{\epsilon})^{1/2-\kappa}})$ is replaced by $1$.

If we look at the estimates of Theorems \ref{cvforte} and \ref{cvfaible} at time $n=n_0$, the factor $\frac{1}{n-1}$ is of size $\Delta t= \text{O}(\Delta t^{1-r})$.

\subsection{Some comments on the convergence results}\label{commsect}

The proofs rely on the following decomposition, which explains the origin of the different error terms:
\begin{equation}\label{thedecompo}
X^{\epsilon}(n\Delta t)-X_n=X^{\epsilon}(n\Delta t)-\overline{X}(n\Delta t)+\overline{X}(n\Delta t)-\overline{X}_n+\overline{X}_n-X_n.
\end{equation}
The numerical process $(\overline{X}_n)$ is defined in \eqref{eqmoynum} below: it is the solution of the macrosolver with a known $\overline{F}$, while $(X_n)$ is solution of the macrosolver using $\tilde{F}_n$. The continuous processes $\overline{X}$ and $X^\epsilon$ are defined at the beginning of Section \ref{sectdescri}.


The first term is bounded thanks to the averaging result, using strong and weak order of convergence results - see the article \cite{CEB1}. The second term is the error in a deterministic numerical scheme, for which convergence results are classical; we recall the estimate in Proposition \ref{propoerrdet}. The third term is the difference between the two numerical approximations, and the main task is the control of this part: we show that an extension of the averaging effect holds at the discrete time level, where $\overline{X}_n$ plays the role of an averaged process for $X_n$.

When we look at the Theorems \ref{cvforte} and \ref{cvfaible}, we first remark that we obtain the same kind of bounds as in the finite dimensional case of \cite{E-L-V}. However we notice some differences; they are due both to the infinite dimensional setting and to different proofs.

First, the weak order of the Euler scheme for SPDEs is only $1/2$ - see \cite{deb} - while it is $1$ for SDEs; in the estimate, the strong order never appears, and this is one of the main theoretical advantages of the method. For completeness, we recall that the strong order is $1/4$ in the SPDE case - see \cite{prin} - and $1/2$ for SDEs. In fact, at least in the strictly dissipative case, we are comparing the invariant measures of the continuous fast equation with the invariant measure of its numerical approximation - see Theorem \ref{weakdistance} and Corollary \ref{weakcor}. In the weakly dissipative case, we use the weak error estimates where the constants do not depend on the final time.




We need some regularity on the initial condition $x$ in Theorem \ref{cvfaible}, since we apply the weak order averaging Theorem of \cite{CEB1}; nonetheless the parameter $\theta>0$ does not appear in the different orders of convergence, but only in the constant $C$.

In Theorem \ref{cvfaible}, we do not exactly obtain the same estimate as in \cite{E-L-V}, since $M$ does not appear: as a consequence, the estimate is not improved when $M$ is increasing. Since we are looking at some weak error estimate, this can seem natural; first $M$ appeared in the estimate of \cite{E-L-V} only because of some strong error estimate used in the proof, while here to prove Theorem \ref{cvfaible} we always think with a weak error state of mind. Second, the corresponding term in the weak estimate of \cite{E-L-V} is easily controlled by $\Delta t$, which already appears in another part of the error since we use a first order scheme for the macrosolver - and not a general scheme - so that there is no need for getting a better estimate.

The proofs of Theorems \ref{cvforte} and \ref{cvfaible} are inspired by the ones in \cite{E-L-V}, but are different. The strong error is analysed in a global way, as in the proof of the strong order Theorem in \cite{CEB1}. Moreover we do not need a counterpart of Lemma $2.6$ in \cite{E-L-V}, and we thus think that our method is more natural. For the control of the weak error, we introduce a new appropriate auxiliary function.

As explained in the introduction, we present here some simplified results and proofs. In the stricly dissipative case, we could go further in the analysis of the error, thanks to the initialization procedure of the scheme, and to additional exponential convergence results. But after a very technical work, due to the regularity assumptions made on the nonlinear coefficients $F$ and $G$ we would only obtain a factor $\Delta t^{1-r-\eta}$ instead of $\Delta t$ in the finite dimensional case - where is any $r>0$, which can be chosen as small as necessary, and $\eta$ is defined in Assumption \ref{hypF} and depends on the regularity properties of the nonlinearities $F$ and $G$. Nevertheless, the cost of the scheme remains of the same order even if we do not use those better estimates.



We consider that $\epsilon$ is fixed, and then the numerical method depends on several parameters $N,M,n_T,\delta t,\Delta t$; we would like to give some explicit choices of the parameters for which we have a simple bound, and for which a direct scheme, where the stiffness of the system \eqref{eqgen} is not treated, is less efficient, for the cost defined by
$$\text{cost}=\frac{Mm_0}{\Delta t}.$$
This is the total number of microtime steps for $T=1$.

We take $n=n_0$: then $\frac{1}{n}\leq C\Delta t$.

The time scale separation parameter $\epsilon$ is supposed to be very small, while we fix some tolerance $tol$ for the error in the numerical scheme: more precisely, we want the error to be
\begin{equation}\label{errexpected}
\begin{gathered}
\E|X^\epsilon(n\Delta t)-X_n|\leq C(\epsilon^{1/2-r}+tol) \text{ (strong error)}\\
|\E\Phi(X^\epsilon(n\Delta t))-\E\Phi(X_n)|\leq C(\epsilon^{1-r}+tol) \text{ (weak error)},
\end{gathered}
\end{equation}
where $r$ satisfies either $r<r_s$ - strong error - or $r<r_w$ - weak error - if the parameters $r_s$ and $r_w$ are defined by the conditions $\epsilon^{1/2-r_s}=\text{o}(tol)$ (strong error) and $\epsilon^{1-r_w}=\text{o}(tol)$ (weak error): the numerical error is dominant with respect to the averaging error. We want to show that we can choose the parameters of the scheme such that each term of the estimate of the Theorem \ref{cvforte} are of size $tol$, except the first one, and such that the cost of the scheme is lower than the cost of a direct scheme.

\begin{itemize}
\item We first notice that the choice for $\Delta t$ and $\frac{\delta t}{\epsilon}$ is easy:
\begin{equation}\label{choiceuniversal}
\begin{gathered}
\Delta t\approx tol^{\frac{1}{1-r}}\\
\frac{\delta t}{\epsilon}\approx tol^{\frac{1}{1/2-\kappa}},
\end{gathered}
\end{equation}
where $r<r_s$ or $r<r_w$, and $\kappa<1/2$ are fixed. The choice of other parameters $N,M,n_T$ changes when we look at the strong or at the weak error.

\item We first focus on the strong convergence case. We consider the case when either $M=1$ or $N=1$. We obtain

\begin{center}\label{choicestrong}
\begin{tabular}{|c|c|c|}
    \hline
     & \text{parameters} & \text{cost}\\
    \hline
    $M=1$ & $N\approx tol^{-2+\frac{1}{1-r}-\frac{1}{1/2-\kappa}}$ & \\
    \cline{2-3} & $n_T\approx tol^{-\frac{1}{1/2-\kappa}}\log(tol^{-1}) $ & $tol^{-2-\frac{1}{1/2-\kappa}+\frac{1}{1-r}}$\\
    \hline
    $N=1$ & $M\approx tol^{\frac{1}{1-r}-2}$ &\\
    \cline{2-3} & $n_T\approx \log(tol^{-1})tol^{-\frac{1}{1/2-\kappa}}$ & $tol^{-2-\frac{1}{1/2-\kappa}+\frac{1}{1-r}}\log(tol^{-1})$\\
    \hline
\end{tabular}
\end{center}

Unlike in the finite dimensional situation, the bound depends on a factor $\Delta t^{1-\eta-r}$ - where $\eta$ is linked to the regularity of the nonlinear coefficients $F$ and $G$. The charge of this difference could be taken by either $N$ or $n_T$; the advantage of choosing $n_T$ is that the difference only involves a logarithmic factor on the final cost. As a consequence, we can make no difference between the choices $N=1$ and $M=1$.


We can now make a comparison with the cost coming from the use of a direct scheme: since the strong order of Euler scheme for SPDEs is $1/4-\kappa$, the error can be bounded by $C(\frac{\delta t}{\epsilon})^{1/4-\kappa}$ for some constant $C$. Therefore to have a bound of size $tol$, the time step size $\delta t$ much satisfy $\frac{\delta}{\epsilon}\approx tol^{\frac{1}{1/4-\kappa}}$. This leads to a cost
$$\frac{1}{\delta t}\approx \epsilon^{-1}tol^{-\frac{1}{1/4-\kappa}};$$
we conclude that in this situation the HMM numerical method is better, since the ratio of the cost tends to $0$ when $\epsilon\rightarrow 0$ - under the condition $r<r_s$.


\item We now focus on the weak error estimate. Here $M$ plays no role, so that $M=1$ is the good choice! The time steps $\Delta t$ and $\delta t$ are still given by \eqref{choiceuniversal}. It remains to look for parameters $N$ and $n_T$ such that
$$\frac{1}{N\frac{\delta t}{\epsilon}+1}e^{-cn_T\frac{\delta t}{\epsilon}}\approx \Delta t.$$

Once again, we need to choose $n_T$ and $N$ such that $N\frac{\delta t}{\epsilon}$ or $n_T\frac{\delta t}{\epsilon}$ is large. Since exponential decrease is faster than polynomial decrease, the best choice is $e^{-cn_T\frac{\delta t}{\epsilon}}\approx \Delta t$, and therefore $n_T\approx tol^{-\frac{1}{1/2-\kappa}}\log(tol^{-1})$, while $N=1$.

We then see that $\frac{1}{((n_T-1)\frac{\delta t}{\epsilon})^{1/2-\kappa}}\leq C$, so that the additional factor appearing when only weak dissipativity is satisfied plays no role.


The corresponding cost of the scheme is then of order $\frac{m_0}{\Delta t}\approx \log(tol^{-1})tol^{-\frac{1}{1-r}-\frac{1}{1/2-\kappa}}$.

We can again compare this cost with the cost coming from the use of a direct scheme: the weak order of Euler scheme for SPDEs is $1/2-\kappa$; to obtain the second estimate of \eqref{errexpected}, 
we need a cost
$$\frac{1}{\delta t}\approx \epsilon^{-1}tol^{-\frac{1}{1/2-\kappa}},$$
and again the HMM scheme is better than the direct scheme for such a range of parameters.


\end{itemize}


\section{Preliminary results}
\subsection{Known results about the fast equation and the averaged equation}\label{known}


In this section, we just recall without proof the main results on the fast equation with frozen slow component and on the averaged equation, defined below. Proofs can be found in \cite{Ce-F} for the strict dissipative case, and the extension to the weakly dissipative situation relies on arguments explained below.

If $x\in H$, we define an equation on the fast variable where the slow variable is fixed and equal to $x$:
\begin{equation}\label{eqfig}
\begin{gathered}
dY_x(t,y)=(BY_x(t,y)+G(x,Y_x(t,y)))dt+dW(t),\\
Y_x(0,y)=y.
\end{gathered}
\end{equation}
This equation admits a unique mild solution, defined on $[0,+\infty[$.

Since $Y^\epsilon$ is involved at time $t>0$, heuristically we need to analyse the properties of $Y_x(\frac{t}{\epsilon},y)$, with $\epsilon\rightarrow 0$, and by a change of time we need to understand the asymptotic behaviour of $Y_x(.,y)$ when time goes to infinity.

Under the strict dissipativity Assumption \ref{strictdiss}, we obtain a contractivity of trajectories issued from different initial conditions and driven by the same noise:
\begin{propo}
With \eqref{hypLg}, for any $t\geq0$, $x,y_1,y_2\in H$  we have
$$|Y_{x}(t,y_1)-Y_{x}(t,y_2)|_H\leq e^{-\frac{(\mu-L_g)}{2}t}|y_1-y_2|_H.$$
\end{propo}
Under the weak dissipativity Assumption \ref{weakdiss}, we obtain such an exponential convergence result for the laws instead of trajectories. The proof of this result is not staightforward, and can be found in \cite{deb-hu-tess} - see also \cite{CEB2} for further references.
\begin{propo}\label{propoexpy1y2}
With \eqref{hypdiss}, there exist $c>0$, $C>0$ such that for any bounded test function $\phi$, any $t\geq 0$ and any $y_1,y_2\in H$
\begin{equation}\label{cvexpy1y2}
|\E\phi(Y(t,y_1))-\E\phi(Y(t,y_2))|\leq C\|\phi\|_{\infty}(1+|y_1|^2+|y_2|^2)e^{-ct}.
\end{equation}
\end{propo}

As a consequence, we can show that there exists a unique invariant probability measure associated with $Y_x$, and that the convergence to equilibrium is exponentially fast.

First, let $\nu=\mathcal{N}(0,(-B)^{-1}/2)$ be the centered Gaussian probability measure on $H$ with the covariance operator $(-B)^{-1}/2$ - which is positive and trace-class, thanks to Assumption \ref{hypAB}.

Then $\mu^x$ defined by
\begin{equation}\label{definvmeas}
\mu^x(dy)=\frac{1}{Z(x)}e^{2U(x,y)}\nu(dy),
\end{equation}
where $Z(x)\in]0,+\infty[$ is a normalization constant, is the unique probability invariant measure associated to $Y_x$.
This expression comes from the gradient structure of equation \eqref{eqfig}, given in Assumption \ref{hypG}.

Second, under any of the dissipativity assumptions, the convergence to equilibrium is exponential in the following sense:
\begin{propo}\label{convexp}
If we assume \eqref{hypLg} or \eqref{hypdiss}, there exist constants $C,c>0$ such that for any bounded function $\phi:H\rightarrow \R$ or $\phi:H\rightarrow H$, $t\geq 0$ and $x,y\in H$ we have
$$|\E\phi(Y_x(t,y))-\int_{H}\phi(z)\mu^x(dz)|\leq C\|\phi\|_{\infty}(1+|y|_{H}^{2})e^{-ct}.$$
\end{propo}


Now we define the averaged equation. First we define the averaged nonlinear coefficient $\overline{F}$:
\begin{defi}For any $x\in H$,
\begin{equation}\label{deffbar}
\overline{F}(x)=\int_{H}F(x,y)\mu^x(dy).
\end{equation}
\end{defi}

Using Assumptions \ref{hypF}, \ref{hypG} and the expression of $\mu^x$, it is easy to prove that $\overline{F}$ is bounded and Lipschitz continuous.



Under Assumptions \ref{hypF}, \ref{hypG}, we can easily prove the following properties on $\overline{F}$:
\begin{propo}\label{FbarLip}
There exists $0\leq \eta<1$ and a constant $C$ such that the following directional derivatives of $\overline{F}$ are well-defined and controlled:
\begin{itemize}
\item For any $x\in H$, $h\in H$, $|D\overline{F}(x).h|\leq C|h|_{H}$.
\item For any $x\in H$, $h\in H$, $k\in D(-A)^\eta$, $|D^2\overline{F}(x).(h,k)|\leq C|h|_{H}|k|_{(-A)^\eta}$.
\item For any $x\in H$, $h,k\in H$, $|(-A)^{-\eta}D^{2}\overline{F}(x).(h,k)|\leq C|h|_H|k|_H$.
\end{itemize}
Moreover, $\overline{F}$ is bounded and Lipschitz continuous.
\end{propo}

The last estimate above is a consequence of Assumptions \ref{hypF2} and \ref{hypG2}, and of the following fact:
we have almost surely $W^B(t)\in L^\infty(0,1)$ for any $t\geq 0$, and
\begin{equation}\label{momentLinfty}
\int_{H}|z|_{L^\infty(0,1)}\nu(dz)<+\infty.
\end{equation}
Moreover if $\eta>1/4$ - which is the right condition in the case of linear Laplace operators and of nonlinear Nemytskii operators - we have $\int_{H}|z|_{(-A)^\eta}\nu(dz)=+\infty$: we thus need the restrictive condition in Assumption \ref{hypG2}.


\begin{rem}
Even when $F$ and $G$ are Nemytskii operators, $\overline{F}$ is not such an operator in general.
\end{rem}
Then the averaged equation (see \eqref{eqmoy} in the introduction) can be defined:
$$\frac{d\overline{X}(t)}{dt}=A\overline{X}(t)+\overline{F}(\overline{X}(t)),$$
with initial condition $\overline{X}(0)=x\in H$.
For any $T>0$, this deterministic equation admits a unique mild solution $\overline{X}\in\mathcal{C}([0,T],H)$.

\subsection{Estimates on the numerical solutions}\label{estimnums}

We can give uniform estimates on $X_n$ and $Y_{n,m,j}$, defined by \ref{schémalent} and \ref{schémarapide}:
\begin{lemme}\label{lentenum}
There exists $C>0$ such that we have $\PP$-almost surely
$$|X_n|\leq C(1+|x|),$$
for any $0\leq n\leq n_0$.
\end{lemme}

\underline{Proof}
The linear operator $S_{\Delta t}$ satifies $|S_{\Delta t}|_{\mathcal{L}(H)}\leq \frac{1}{1+\lambda \Delta t}$; moreover $F$ is bounded, so that by \eqref{ftilde} we almost surely have for any $n\geq 0$ $|\tilde{F}_n|\leq \|F\|_{\infty}$. The end of the proof is then straightforward; we also notice that $C$ does not depend on the final time $T$.
\qed

\begin{lemme}\label{rapidenum}
There exists $C>0$ - which does not depend on $T>0$, on $N$, on $n_T$ or on $M$ - such that for any $\Delta t>0$, $\tau=\frac{\delta t}{\epsilon}>0$, $0\leq n\leq n_0$, $0\leq m\leq m_0$ and $1\leq j\leq M$, we have
$$\E|Y_{n,m,j}|^2\leq C(|y|^2+1).$$
\end{lemme}

\underline{Proof}
We introduce $\omega_{n,m,j}$ defined by the fast numerical scheme with no nonlinear coefficient - see \eqref{schémarapide} - with the notation $\tau=\frac{\delta t}{\epsilon}$: for any $0\leq j\leq M$, $0\leq n\leq n_0$ and $0\leq m\leq m_0$
\begin{equation}\label{defnumsto}
w_{n,m+1,j}=R_{\tau}w_{n,m,j}+\sqrt{\tau}R_{\tau}\zeta_{n,m+1,j},
\end{equation}
with the initial condition $w_{n,0,j}=0$.
$w_{n,m,j}$ is the numerical approximation - using the microsolver scheme with a step size $\tau$ - of the process defined by
\begin{gather*}
d\textbf{Z}^{(n,j)}(t)=B\textbf{Z}^{(n,j)}(t)dt+d\tilde{W}^{(n,j)}(t)\\
\textbf{Z}^{(n,j)}(0)=0,
\end{gather*}
where $\tilde{W}_{t}^{(n,j)}=\frac{1}{\sqrt{\epsilon}}W_{\epsilon t}^{(n,j)}$. Notice that the $(\tilde{W}^{(n,j)})$ are independent cylindrical Wiener processes, and that
$$\zeta_{n,m+1,j}=\frac{\tilde{W}_{(m+1)\tau}^{(n,j)}-\tilde{W}_{m\tau}^{(n,j)}}{\sqrt{\tau}}.$$
The $\textbf{Z}^{(n,j)}$ are independent realizations of $Z(t)=W^B(t)=\int_{0}^{t}e^{(t-s)B}dW(s)$. Using Theorem $3.2$ of \cite{prin}, giving the strong order $1/4$ for the microsolver - when the initial condition is $0$, with no nonlinear coefficient, with a constant diffusion term and under the assumptions made here - we get the following estimate: for any $0<\kappa<1/2$, there exists $C_{\kappa}>0$, such that for any $\tau>0$, $0\leq n\leq n_0$, $1\leq j\leq M$ and $0\leq m\leq m_0$
\begin{equation}\label{estimnoisenum}
\E|\omega_{n,m,j}-\textbf{Z}^{(n,j)}(m\tau)|^{2}\leq C_{\kappa}\tau^{1/2-\kappa}.
\end{equation}
Thanks to \eqref{estimWB} and \eqref{estimnoisenum}, for any $\tau_0>0$, there exists a constant $C(\tau_0)$ such that for any $0\leq \tau<\tau_0$, $0\leq n\leq n_0$, $1\leq j\leq M$ and $0\leq m\leq m_0$ we have $$\E|\omega_{n,m,j}|^2\leq C(\tau_0).$$
Now for any $0\leq m\leq m_0$ we define $D_{n,m,j}=Y_{n,m,j}-w_{n,m,j}$; it is enough to control $|D_{n,m,j}|^2$, $\PP$-almost surely. By \eqref{schémarapide}, we have the following expression: for any $0\leq m\leq m_0$
$$D_{n,m+1,j}=R_{\tau}D_{n,m,j}+\tau R_{\tau}G(X_n,Y_{n,m,j}).$$
Since $G$ is bounded, and using the inequality $|R_\tau|_{\mathcal{L}(H)}\leq \frac{1}{1+\mu\tau}$, we get
$$|D_{n,m+1,j}|\leq\frac{1}{1+\mu\tau}|D_{n,m,j}|+C\tau.$$
Therefore we have for any $0\leq m\leq m_0$
$$(1+\mu\tau)^m|D_{n,m,j}|\leq |D_{n,0,j}|+C[(1+\mu\tau)^m-1].$$
But $D_{n,0,j}=D_{n-1,m_0,j}$. So we get
$$(1+\mu\tau)^{m_0}|D_{n+1,0,j}|\leq |D_{n,0,j}|+C[(1+\mu\tau)^{m_0}-1].$$
Therefore
$$(1+\mu\tau)^{nm_0}|Z_{n,0,j}|\leq |D_{0,0,j}|+C(1+\mu\tau)^{nm_0}=|y|+C(1+\mu\tau)^{nm_0}.$$
As a consequence, we get for any $0\leq n\leq n_0$ and any $0\leq m\leq m_0$, $|D_{n,0,j}|\leq C(1+|y|)$, and then $|Z_{n,m,j}|\leq 2C(1+|y|)$.
\qed

\subsection{Asymptotic behaviour of the ``fast'' numerical scheme}\label{invarfast}



At the continuous time level, the averaging principle proved in \cite{CEB1} comes from the asymptotic behaviour of the fast equation with frozen slow component \eqref{eqfig}, as described in Section \ref{known}. The underlying idea of the HMM method in our setting is to prove a similar averaging effect at the discrete time level: we therefore study the asymptotic behaviour of the fast numerical scheme which defines the microsolver, with frozen slow component - in other words we are looking at the evolution of the microsolver during one fixed macrotime step.

In Section \ref{known}, we have seen that under the weak dissipativity Assumption \ref{weakdiss} the fast equation with frozen slow component admits a unique invariant probability measure $\mu^x$. At the discrete time level, this Assumption only yields the existence of invariant laws; to get a unique invariant law $\mu^{x,\tau}$, we need strict dissipativity \ref{hypLg}, and we obtain:
\begin{theo}\label{unicinv}
Under Assumption \ref{strictdiss}, for any $\tau>0$ and any $x\in H$; the numerical scheme \ref{schémarapide2} admits a unique ergodic invariant probability measure $\mu^{x,\tau}$.
Moreover, we have convergence to equilibrium in the following sense: for any $\tau_0>0$, there exist $c>0$ and $C>0$ such that for any $0<\tau\leq \tau_0$, $x\in H$, $y\in H$, any Lipschitz continuous function $\phi$ from $H$ to $R$, and $m\geq 0$, we have
\begin{equation}\label{estimvitesse}
|\E(\phi(Y_{m}^{x}(y)))-\int_{H}\phi(z)\mu^{x,\tau}(dz)|\leq C(1+|y|)[\phi]_{\text{Lip}}e^{-cm\tau}.
\end{equation}
\end{theo}

We recall the notation $\tau=\frac{\delta t}{\epsilon}$ for the effective time step; the noise is defined with a cylindrical Wiener process $\tilde{W}$: $\tilde{\zeta}_{m+1}=\frac{\tilde{W}_{(m+1)\tau}-\tilde{W}_{m\tau}}{\sqrt{\tau}}$. If we fix the slow component $x\in H$, we define
\begin{equation}\label{schémarapide2}
Y_{m+1}^{x}(y)=R_{\tau}Y_{m}^{x}(y)+\tau R_{\tau}G(x,Y_{m}^{x}(y))+\sqrt{\tau}R_{\tau}\tilde{\zeta}_{m+1},
\end{equation}
with the initial condition $Y_{0}^{x}(y)=y$.

\subsubsection{Existence of an invariant law}

With Equation \eqref{schémarapide2}, we associate the transition semi-group $(P_{m}^{x,\tau})$: if $\phi$ is a bounded measurable function from $H$ to $\R$, if $y\in H$ and $m\geq 0$
\begin{equation}\label{defPm}
P_{m}^{x,\tau}\phi(y)=\E[\phi(Y_{m}^{x}(y))].
\end{equation}
We also denote by $\nu_{m,y}^{x,\tau}$ the law of $Y_{m}^{x,\tau}(y)$: then
$$P_{m}^{x,\tau}\phi(y)=\E[\phi(Y_{m}^{x,\tau}(y))]=\int_{H}\phi(z)\nu_{m,y}^{x,\tau}(dz).$$
We notice that the semi-group $(P_m)$ satisfies the Feller property: if $\phi$ is bounded and continuous, then $P_m\phi\in$ is bounded and continuous.

The required tightness property is a consequence of the following estimate, which can be proved thanks to regularization properties of the semi-group $(R_{\tau}^{m})_m$: for any $0<\gamma<1/4$, $\tau>0$, there exists $C(\gamma,\tau)>0$ such that for any $m\geq 1$ and $\tau\leq \tau_0$
$$\E|Y_{m}(y)|_{(-B)^\gamma}^{2}\leq C(\gamma,\tau).$$
Moreover if $0<\gamma<1/4$, the embedding of $D(-B)^\gamma$ in $H$ is compact.

\subsubsection{Uniqueness under strict dissipativity}

The key estimate to prove uniqueness is the following contractivity property, which holds thanks to Assumption \ref{strictdiss}:
\begin{propo}\label{theoestim}
For any $\tau_0>0$, there exists $c>0$ such that for any $0<\tau\leq \tau_0$, $m\geq 0$, $y_1,y_2\in H$, $x\in H$ we have $\PP$-almost surely
\begin{equation}\label{theestim}
|Y_{m}^{x}(y_1)-Y_{m}^{x}(y_2)|\leq e^{-cm\tau}|y_1-y_2|.
\end{equation}
\end{propo}
\underline{Proof}
If we define $r_m=Y_m(y_1)-Y_m(y_2)$, then we have the equation
\begin{gather*}
r_{m+1}=r_{m}+\tau BR_{m+1}+\tau (G(x,Y_{m}(y_1))-G(x,Y_{m}(y_2))),\\
r_{0}=y_1-y_2.
\end{gather*}
If we take the scalar product in $H$ of this equation with $r_{m+1}$, we get
\begin{align*}
|r_{m+1}|^2-<r_{m},r_{m+1}>&=\tau<Br_{m+1},r_{m+1}>\\
&+\tau <G(x,Y_{m}(y_1))-G(x,Y_{m}(y_2)),r_{m+1}>;
\end{align*}
The left hand-side is equal to $\frac{1}{2}(|r_{m+1}|^2-|r_{m}|^2)+\frac{1}{2}|r_{m+1}-r_{m}|^2,$ and we get
\begin{align*}
\frac{1}{2}(|r_{m+1}|^2-|r_{m}|^2)&\leq -\tau|(-B)^{1/2}r_{m+1}|^2+\tau L_g|r_{m}||r_{m+1}|\\
&\leq -\mu\tau|r_{m+1}|^2+\frac{1}{2}\tau L_g(|r_{m+1}|^2+|r_{m}|^2).
\end{align*}
Therefore we have $(1+\tau(2\mu-L_g))|r_{m+1}|^2\leq(1+\tau L_g)|r_{m}|^2.$
We remark that \ref{hypLg} implies that for any $\tau_0>0$, there exists $c>0$ such that if $\tau\leq \tau_0$ we have $\rho=\frac{1+\tau L_g}{1+\tau(2\mu-L_g)}\leq e^{-2c\tau}$; therefore
$$|r_{m}|^2\leq \rho^m|y_1-y_2|^2\leq e^{-2cm\tau}|y_1-y_2|^2.$$
\qed

As a consequence, there exists a unique ergodic invariant probability measure $\mu^{x,\tau}$, which is strongly mixing. Moreover one can prove that there exists $C>0$ such that for any $\tau>0$ and any $x\in H$ $\int_{H}|y|\mu^{x,\tau}(dy)\leq C$; we therefore get
\begin{align*}
|\E(\phi(Y_{m}^{x}(y)))-\int_{H}\phi(z)\mu^{x,\tau}(dz)|&=|\E(\phi(Y_{m}^{x}(y)))-\int_{H}\E\phi(Y_{m}^{x}(z))\mu^{x,\tau}(dz)|\\
&\leq\int_{H}\E|\phi(Y_{m}^{x}(y))-\phi(Y_{m}^{x}(z)|\mu^{x,\tau}(dz)\\
&\leq[\phi]_{\text{Lip}}\int_{H}e^{-cm\tau}|y-z|\mu^{x,\tau}(dz)\\
&\leq C[\phi]_{\text{Lip}}e^{-cm\tau}.
\end{align*}

\subsubsection{Approximation of the invariant law $\mu^x$ by the fast numerical scheme}


We recall that $\mu^{x}$ denotes the invariant law of the continuous time fast equation with frozen slow component \eqref{eqfig}. Thanks to the fast numerical scheme, we can get an approximation result, which is proved in \cite{CEB2}: with test functions of class $\mathcal{C}_{b}^{2}$, we can control the weak error for any time with a convergence of order $1/2$ with respect to the time step $\tau$. Moreover the estimate is easily seen to be independent from the slow component $x$.

We define for any $\Phi$ of class $\mathcal{C}_{b}^{2}$
$$\|\Phi\|_{(2)}=\sup_{y\in H}|\Phi(y)|+\sup_{y\in H,h\in H,|h|=1}|D_y\Phi(y).h|+\sup_{y\in H,h,k\in H,|h|=|k|=1}|D_{yy}^{2}\Phi(y).(h,k)|.$$

\begin{theo}\label{weakdistance}
With the dissipativity condition \eqref{hypdiss}, for any $0<\kappa<1/2$, for any $\tau_0>0$, there exists $C,c>0$ such that for any $\Phi$ of class $\mathcal{C}_{b}^{2}$, for any $x,y\in H$, for any $0<\tau\leq \tau_0$ and any integer $2\leq m< +\infty$
$$|\E[\Phi(Y_x(m\tau,y))]-\E[\Phi(Y_{m}^{x}(y))]|\leq C\|\Phi\|_{(2)}(1+|y|^3)(((m-1)\tau)^{-1/2+\kappa}+1)\tau^{1/2-\kappa}.$$
\end{theo}

As explained in Section \ref{invarfast}, in general the fast numerical scheme does not admit a unique invariant probability measure, but when the strict dissipativity Assumption \ref{strictdiss} is satisfied, it admits a unique invariant law $\mu^{x,\tau}$.

\begin{cor}\label{weakcor}
Under the assumptions of Theorem \ref{weakdistance}:

\em(i) if only \eqref{hypdiss} is satisfied, we have for any $m\geq 2$
$$|\int_{H}\Phi d\mu^x-\E[\Phi(Y_{m}^{x})]|\leq C\|\Phi\|_{(2)}(1+|y|^3)(((m-1)\tau)^{-1/2+\kappa}+1)\tau^{1/2-\kappa}+CN(\Phi)(1+|y|^2)e^{-cm\tau}.$$

\em(ii) if moreover \eqref{hypLg} is satisfied,
$$|\int_{H}\Phi d\mu^{x}-\int_{H}\Phi d\mu^{x,\tau}|\leq C\|\Phi\|_{(2)}\tau^{1/2-\kappa}.$$
\end{cor}

We recall that in the case of Euler scheme for SDEs this kind of results holds with the order of convergence $1$.


\subsection{Error in the deterministic scheme \eqref{eqmoynum}}\label{errdet}


We define a scheme based on the macrosolver, for theoretical purpose, in the situation when $\overline{F}$ is known:
\begin{equation}\label{eqmoynum}
\begin{gathered}
\overline{X}_{n+1}=S_{\Delta t}\overline{X}_n+\Delta t S_{\Delta t}\overline{F}(\overline{X}_n)\\
\overline{X}_{0}=x.
\end{gathered}
\end{equation}

We can look at the error between $\overline{X}_n$, defined by \eqref{eqmoynum}, and $\overline{X}(n\Delta t)$, defined by \eqref{eqmoy}. Here quantities are deterministic, and the following result is classical - see \cite{le_roux}, \cite{crou-tho}, or the details of the proofs in \cite{prin}: 

\begin{propo}\label{propoerrdet}For any $0<r<1$, $\Delta t_0>0$ and $T>0$, there exists $C>0$, such that for any $0<\Delta t\leq \Delta t_0$ and $1\leq n\leq \lfloor \frac{T}{\Delta t}\rfloor$
$$|\overline{X}_n-\overline{X}(n\Delta t)|\leq \frac{C}{n}+C(1+|x|)\Delta t^{1-r}.$$
\end{propo}

\section{Proof of the strong convergence Theorem \ref{cvforte}}\label{sectfort}

Thanks to the following remark on the construction of the scheme, we can consider that $Y_{n,m,j}$ corresponds to a process evaluated at time $m+nm_0$. The idea is to build noise processes $\zeta^{(j)}$ by concatenation of the $\zeta^{(n,j)}$ for the different $n$.
\begin{rem}\label{remarknoise}
To compute $Y_{n,m,j}$ with the microsolver, the total number of microtime steps used is $nm_0+m$. It is then natural to use only one noise process $\zeta_{.}^{(j)}$, for each $1\leq j\leq M$, and to make an evaluation involving time $nm_0+m$. More precisely, we can use
$$\zeta_{k+1}^{(j)}=\frac{W_{(k+1)\delta t}^{(j)}-W_{k\delta t}^{(j)}}{\sqrt{\delta t}},$$
with $(W^{(j)})_{1\leq j\leq M}$ being $M$ independent cylindrical Wiener processes on $H$.

Then one can define for $0\leq t\leq m_0\delta t$ and for $0\leq m\leq m_0-1$
\begin{gather*}
W^{(n,j)}(t)=W^{(j)}(t+nm_0\delta t)-W^{(j)}(nm_0\delta t)\\
\zeta_{n,m+1,j}=\zeta_{nm_0+m+1}^{(j)}.
\end{gather*}
\end{rem}

We also introduce the following notation: $\E_n$ denotes conditional expectation with respect to the $\sigma$-field
\begin{equation}\label{defsigmafield}
\mathcal{G}_n=\sigma(\zeta_{k,m,j},0\leq k\leq n-1,1\leq m\leq m_0,1\leq j\leq M)=\sigma(\zeta_{m}^{(j)},1\leq m\leq nm_0).
\end{equation}

We notice that $X_n$ is $\mathcal{G}_n$-measurable, but that $\tilde{F}_n$ is not.

The final time $T$ is fixed and we recall the notation $n_0=\lfloor\frac{T}{\Delta t}\rfloor$.

To simplify notations, we do not always precise the range of summation in the expressions below: the indices $j,j_1,j_2$ belong to $\left\{1,\ldots,M\right\}$, and $m,m_1,m_2$ belong to $\left\{n_T,\ldots,n_T+N-1=m_0\right\}$.

We recall that according to the decomposition of the error \eqref{thedecompo}, we have to control
\begin{equation}\label{thedecompostrong}
\begin{aligned}
\E|X^\epsilon(n\Delta t)-X_n|&\leq \E|X^{\epsilon}(n\Delta t)-\overline{X}(n\Delta t)|\\
&\hspace{10 pt}+|\overline{X}(n\Delta t)-\overline{X}_n|\\
&\hspace{10 pt}+\E|\overline{X}_n-X_n|.
\end{aligned}
\end{equation}
The first part is controlled thanks to the strong order Theorem of \cite{CEB1}: for any $0<r<1/2$, we have for any $0\leq n\leq n_0$
$$\E|X_{n\Delta t}^{\epsilon}-\overline{X}_{n\Delta t}|\leq C\epsilon^{1/2-r}.$$
The second part is deterministic and is controlled thanks to Proposition \ref{propoerrdet}:
$$|\overline{X}_n-\overline{X}(n\Delta t)|\leq \frac{C}{n}+C(1+|x|)\Delta t^{1-r},$$
where $C$ depends on $T,r,x,y$.

It remains to focus on the third part $e_n=\overline{X}_n-X_n$. Instead of analyzing the local error like in \cite{E-L-V}, we adopt a global point of view, and we follow the idea of the proof of Theorem $1.1$ in \cite{CEB1}: for any $0\leq n\leq n_0$
\begin{equation}\label{globalerr}
X_n-\overline{X}_n=S_{\Delta t}^{n}x+\Delta t\sum_{k=0}^{n-1}S_{\Delta t}^{n-k}\tilde{F}_k-S_{\Delta t}^{n}x-\Delta t\sum_{k=0}^{n-1}S_{\Delta t}^{n-k}\overline{F}(\overline{X}_k).
\end{equation}
The averaged coefficient $\overline{F}$ is Lipschitz continuous, and $|S_{\Delta t}|_{\mathcal{L}(H)}\leq 1$; moreover we can define the averaged coefficient $\overline{F}^{\tau}$ with respect to the invariant measure $\mu^{x,\tau}$ of the fast numerical scheme - which is unique since we assume strict dissipativity \eqref{hypLg}: for any $x\in H$
\begin{equation}\label{deffbartau}
\overline{F}^{\tau}(x)=\int_{H}F(x,y)\mu^{x,\tau}(dy).
\end{equation}

The error in \eqref{globalerr} can then be decomposed in the following way - the idea of looking at the square of the norm in the second expression is an essential tool of the proof:
\begin{equation}\label{globalerr2}
\begin{aligned}
\E|X_n-\overline{X}_n|&\leq C\Delta t\sum_{k=0}^{n-1}\E|X_k-\overline{X}_k|\\
&\hspace{10 pt}+\left(\E|\Delta t\sum_{k=0}^{n-1}S_{\Delta t}^{n-k}\tilde{F}_k-S_{\Delta t}^{n-k}\overline{F}^\tau(X_k)|^2\right)^{1/2}\\
&\hspace{10 pt}+\Delta t\sum_{k=0}^{n-1}\E|S_{\Delta t}^{n-k}\overline{F}(X_k)-S_{\Delta t}^{n-k}\overline{F}^\tau(X_k)|.
\end{aligned}
\end{equation}
If we can control the two last terms by a certain quantity $Q$, by a discrete Gronwall Lemma we get for any $0\leq k\leq n_0$ $\E|X_n-\overline{X}_n|\leq e^{CT}Q$.

First, the third term in \eqref{globalerr2} is linked to the distance between the invariant measures $\mu^x$ and $\mu^{x,\tau}$ - since we assume strict dissipativity for this strong estimate - which is evaluated thanks to Theorem \ref{weakdistance} and Corollary \ref{weakcor} for test functions of class $\mathcal{C}_{b}^{2}$. Since $|(-A)^\eta S_{\Delta t}^{n-k}|_{\mathcal{L}(H)}\leq \frac{C}{((n-k)\Delta t)^\eta}$, using Assumption \ref{hypF2}, we can apply Corollary \ref{weakcor} to obtain
$$|S_{\Delta t}^{n-k}\overline{F}(X_k)-S_{\Delta t}^{n-k}\overline{F}^\tau(X_k)|\leq \frac{C}{((n-k)\Delta t)^\eta}\tau^{1/2-\kappa},$$
and summing we get
$$\Delta t\sum_{k=0}^{n-1}\E|S_{\Delta t}^{n-k}\overline{F}(X_k)-S_{\Delta t}^{n-k}\overline{F}^\tau(X_k)|\leq C\tau^{1/2-\kappa}.$$

The other term is more complicated; in order to get a precise estimate, we expand the square of the norm of the sum. We can then use some conditional expectations, which allow to use exponential convergence to equilibrium via Theorem \ref{unicinv}. Therefore we obtain the following expansion:
\begin{align*}
\E|\Delta t\sum_{k=0}^{n-1}S_{\Delta t}^{n-k}\tilde{F}_k-S_{\Delta t}^{n-k}&\overline{F}^\tau(X_k)|^2=\Delta t^2\sum_{k=0}^{n-1}\E|S_{\Delta t}^{n-k}(\tilde{F}_k-\overline{F}^\tau(X_k))|^2\\
&+2\Delta t^2\sum_{0\leq k_1<k_2\leq n-1}\E<S_{\Delta t}^{n-k_1}(\tilde{F}_{k_1}-\overline{F}^\tau(X_{k_1})),S_{\Delta t}^{n-k_2}(\tilde{F}_{k_2}-\overline{F}^\tau(X_{k_2}))>\\
&=:\Sigma_1+\Sigma_2.
\end{align*}
{\em (i)} We first treat $\Sigma_1$.

By \eqref{ftilde}, for any $0\leq k\leq n-1$ we have $$\tilde{F}_k=\frac{1}{MN}\sum_{j=1}^{M}\sum_{m=n_T}^{n_T+N-1}F(X_k,Y_{k,m,j});$$
therefore we can see that
\begin{multline*}
\E|S_{\Delta t}^{n-k}(\tilde{F}_k-\overline{F}^\tau(X_k))|^2\\
=\frac{1}{M^2N^2}\sum_{j_1,j_2}\sum_{m_1,m_2}\E\E_k<S_{\Delta t}^{n-k}(F(X_k,Y_{k,m_1,j_1})-\overline{F}^{\tau}(X_k)),S_{\Delta t}^{n-k}(F(X_k,Y_{k,m_2,j_2})-\overline{F}^\tau(X_k))>,
\end{multline*}
with the conditional expectation $\E_k$ with respect to $\mathcal{G}_k$ - see \eqref{defsigmafield}.

When $j_1\neq j_2$, $\zeta_{m}^{(j_1)}$ and $\zeta_{m}^{(j_2)}$ are independent, so that we treat differently the cases $j_1=j_2$ and $j_1\neq j_2$ in the above summation, and we obtain
\begin{multline*}
M^2N^2\E|S_{\Delta t}^{n-k}(\tilde{F}_k-\overline{F}^\tau(X_k))|^2\\
=\sum_{j_1\neq j_2}\E<\sum_{m_1}\E_kS_{\Delta t}^{n-k}(F(X_k,Y_{k,m_1,j_1})-\overline{F}^{\tau}(X_k)),\sum_{m_2}\E_kS_{\Delta t}^{n-k}(F(X_k,Y_{k,m_2,j_2})-\overline{F}^{\tau}(X_k))>\\
+\sum_{j=1}^{M}\sum_{m_1,m_2}\E<S_{\Delta t}^{n-k}(F(X_k,Y_{k,m_1,j})-\overline{F}^{\tau}(X_k)),S_{\Delta t}^{n-k}(F(X_k,Y_{k,m_2,j})-\overline{F}^\tau(X_k))>.
\end{multline*}

For the first part, we can directly use the exponential convergence to equilibrium result of \eqref{estimvitesse} to have a bound with
\begin{align*}
(\frac{1}{N}\sum_{m=n_T}^{m_0}e^{-cm\tau})^2&\leq (\frac{1}{N}\frac{e^{-cn_T\tau}-e^{-c(N+n_T)\tau}}{1-e^{-c\tau}})^2\\
&\leq (Ce^{-cn_T\tau}\frac{1-e^{-cN\tau}}{N\tau})^2\\
&\leq (\frac{Ce^{-cn_T\tau}}{N\tau+1})^2.
\end{align*}
For the second part, we see that we can treat only the case $m_1\leq m_2$, and we introduce the conditional expectation $\E_{k,m_1,j}$ with respect to the $\sigma$-field generated by $\mathcal{G}_k$ and $(\zeta_{km_0+m}^{(j)})_{0\leq m\leq m_1-1}$, when $m_1\leq m_2$. We get a bound with
$$\frac{2}{MN^2}\sum_{n_T\leq m_1\leq m_2\leq m_0}e^{-c(m_2-m_1)\tau}\leq \frac{C}{M(N\tau+1)}.$$
We therefore get
\begin{equation}\label{estimSigma_1}
\Sigma_1\leq C\Delta t\bigl((\frac{e^{-cn_T\tau}}{N\tau+1})^2+\frac{1}{M(N\tau+1)}\bigr).
\end{equation}


{\em (ii)} We now consider $\Sigma_2$, which corresponds to the cross-terms in the expansion of the square of the norm of the quantity $S_{\Delta t}^{n-k}(\tilde{F}_k-\overline{F}^\tau(X_k))$. By the definition of $\tilde{F}_k$, the general term with indices $k_1<k_2$ in $|\Sigma_2|$ is bounded by
\begin{multline*}
|\E<S_{\Delta t}^{n-k_1}(\tilde{F}_{k_1}-\overline{F}^\tau(X_{k_1})),S_{\Delta t}^{n-k_2}(\tilde{F}_{k_2}-\overline{F}^\tau(X_{k_2}))>|\\
\leq\frac{\Delta t^2}{M^2N^2}|\sum_{m_i,j_i}\E<S_{\Delta t}^{n-k_1}(F(X_{k_1},Y_{k_1,m_1,j_1})-\overline{F}^{\tau}(X_{k_1})),S_{\Delta t}^{n-k_2}(F(X_{k_2},Y_{k_2,m_2,j_2})-\overline{F}^{\tau}(X_{k_2}))>|\\
\leq C\frac{\Delta t^2}{MN}\sum_{m=n_T}^{m_0}\E|\E_{k_2}[S_{\Delta t}^{n-k_2}(F(X_{k_2},Y_{k_2,m,j})-\overline{F}^\tau(X_{k_2})]|,
\end{multline*}
using conditional expectation $\E_{k_2}$ and the boundedness of $F$.

Using the exponential convergence result of \eqref{estimvitesse} and Lemma \ref{rapidenum}, we get the bound
$$\E|\E_{k_2}[S_{\Delta t}^{n-k_2}(F(X_{k_2},Y_{k_2,m,j})-\overline{F}^\tau(X_{k_2})]|\leq Ce^{-m\tau},$$
so that the previous quantity is bounded by
$$C\frac{\Delta t^2}{MN}\sum_{m=n_T}^{m_0}\sum_{j=1}^{M}e^{-cm\tau}\leq C\Delta t^2\frac{e^{-cn_T\tau}}{N\tau+1}.$$

Summing on $k_1<k_2$, we can now conclude that
\begin{equation}\label{estimSigma_2}
\Sigma_2\leq C\frac{e^{-cn_T\tau}}{N\tau+1};
\end{equation}
then by \eqref{estimSigma_1} and \eqref{estimSigma_2}
$$
\E|\Delta t\sum_{k=0}^{n-1}S_{\Delta t}^{n-k}\tilde{F}_k-S_{\Delta t}^{n-k}\overline{F}^\tau(X_k)|^2
\leq C\Bigl(\frac{e^{-cn_T\tau}}{N\tau+1}+\Delta t(\frac{e^{-cn_T\tau}}{N\tau+1})^2+\frac{\Delta t}{M(N\tau+1)}\Bigr),
$$
and the result of Theorem \ref{cvforte} now follows from \eqref{globalerr2}.

\section{Proof of the weak convergence Theorem \ref{cvfaible}}\label{sectfaible}


In order to get a better bound for the weak error than for the strong error, we use an auxiliary function which is solution of a Kolmogorov equation.

The proof is valid under the general dissipativity Assumption \ref{weakdiss}.

We divide the proof in two parts: the first one contains the elements of the proof, while the second one is devoted to two technical lemmas.

\subsection{Proof of the Theorem}

According to the decomposition \eqref{thedecompo}, we want to control for any $0\leq n\leq n_0$
\begin{equation}\label{thedecompoweak}
\begin{aligned}
|\E\Phi(X^{\epsilon}(n\Delta t))-\E\Phi(X_n)|&\leq |\E\Phi(X^{\epsilon}(n\Delta t))-\E\Phi(\overline{X}(n\Delta t))|\\
&\hspace{10 pt}+|\Phi(\overline{X}(n\Delta t))-\Phi(\overline{X}_{n})|\\
&\hspace{10 pt}+|\Phi(\overline{X}_{n})-\E\Phi(X_n)|.
\end{aligned}
\end{equation}

Thanks to the averaging Theorem of \cite{CEB1}, which is proved under the weak dissipativity assumption \eqref{hypdiss}, the first term above can be controlled by $C_r\epsilon^{1-r}$, where $C_r$ is a constant - depending on $r,\Phi,x,y,T$, for any $0<r<1$.

For the second term, since we look at the error made by using a deterministic scheme to approximate a deterministic equation, there is no difference between the strong and the weak orders - since the test function $\Phi$ is Lipschitz continuous; so we again use Proposition \ref{propoerrdet}.

For the third term, we see that we have to control some error between two different numerical schemes, in a weak sense. The usual strategy is to decompose this error by means of an auxiliary function satisfying some kind of Kolmogorov equation.

More precisely, we use the deterministic scheme defining $\overline{X}_k$ in order to define for any $0\leq k\leq n$
\begin{equation}\label{fonctionaux}
u_n(k,x)=\Phi(\overline{X}_{n-k}(x)),
\end{equation}
where we explicitly mention dependence of the numerical solution $\overline{X}_k$ on the initial condition $x$.

\begin{rem}
We can easily prove that we have
\begin{gather*}
u_n(n,x)=\Phi(x)\\
u_n(k,x)=u(k+1,S_{\Delta t}x+\Delta t S_{\Delta t}\overline{F}(x)) \text{ for any }k<n.
\end{gather*}
This is the way this function is defined in \cite{E-L-V}.
\end{rem}


We now analyse the error by identifying a telescoping sum:
\begin{equation}\label{decompoweak2}
\begin{aligned}
|\Phi(\overline{X}_n)-\E\Phi(X_n)|&=|u_n(0,x)-\E u_n(n,X_n)|\\
&=|\sum_{k=0}^{n-1}(\E u_n(k,X_k)-\E u_n(k+1,X_{k+1}))|\\
&\leq \sum_{k=0}^{n-1}|\E u_n(k,X_k)-\E u_n(k+1,X_{k+1})|.
\end{aligned}
\end{equation}

According to Lemma \ref{regu_n} below, $u_n$ is of class $\mathcal{C}_{b}^{2}$, and we can control the first and second order derivatives:
\begin{lemme}\label{regu_n}
For any $0< T<+\infty$, there exists $C_{T}>0$ such that for any $0\leq n\leq n_0=\lfloor \frac{T}{\Delta t} \rfloor$ and $0\leq k\leq n$, we have, for any $x\in H$, $h\in H$, $h_1,h_2\in H$:
\begin{gather*}
|D_xu_n(k,x).h|\leq C_{T}|h|,\\
|D_{xx}^{2}u_n(k,x).(h_1,h_2)|\leq C_{T}|h_1||h_2|.
\end{gather*}
Moreover for any $0\leq k\leq n-1$ and any $x\in H$, $h\in H$,
$$|D_xu_n(k,x).h|\leq C_{T}(\Delta t|h|+\frac{|h|_{(-A)^{-\eta}}}{((n-k)\Delta t)^\eta}),$$
where $\eta$ is defined in Assumption \ref{hypF}. 
\end{lemme}

Since the auxiliary functions $u_n$ are linked to the deterministic discrete-time process $(\overline{X}_n)_{n\geq 0}$, the proof does not use stochastic tools.

Moreover the second and the third estimates of this Lemma reveal some smoothing effect in the equation, due to the semi-group $(S_{\Delta t}^{n})_{n\in \N}$. These results are specific to the infinite dimensional framework; if we only use the first estimate above, and a simple one like
$$|D_{xx}^{2}u_n(k,x).(h_1,h_2)|\leq C_{T}|h_1||h_2|_{(-A)^\eta},$$
we would not obtain an optimal convergence result.

To make the proof of Theorem \ref{cvfaible} clearer, we postpone the proof of Lemma \ref{regu_n} in Section \ref{sect_aux_proofs}.

In the general term of \eqref{decompoweak2}, we proceed with a Taylor expansion and we get
\begin{equation}\label{taylorexpansion}
\begin{aligned}
|\E u_n(k,X_k)-&\E u_n(k+1,X_{k+1})|\\
&=|\E u_n(k+1,S_{\Delta t}X_k+\Delta t S_{\Delta t}\overline{F}(X_k))-\E u_n(k+1,S_{\Delta t}X_k+\Delta t S_{\Delta t}\tilde{F}_k)|\\
&\leq \Delta t|\E D_xu_n(k+1,S_{\Delta t}X_k+\Delta t S_{\Delta t}\overline{F}(X_k)).(S_{\Delta t}\overline{F}(X_k)-S_{\Delta t}\tilde{F}_k)|\\
&\hspace{10 pt}+C\Delta t^2\E|\tilde{F}_k-\overline{F}(X_k)|^2.
\end{aligned}
\end{equation}
The second term is controlled by using the second estimate of the previous Lemma \ref{regu_n}. Under the assumption that $F$ is bounded, we therefore get a $\text{O}(\Delta t^2)$ term; when summing over $0\leq k\leq n-1$, we get a $\text{O}(\Delta t)$ term, which is already dominated in the final estimate.

When $k=n-1$, since $u_n(n,.)=\Phi$,
\begin{equation}\label{case_n-1}
|\E u_n(n-1,X_{n-1})-\E u_n(n,X_{n})|\leq C\Delta t.
\end{equation}
In the rest of the proof, we focus on the general case $0\leq k<n-1$.

For the first term, we do not exactly follow the proof of \cite{E-L-V}; we rather define auxiliary functions for $0\leq k\leq n$ in order to keep on looking at a weak error term:
\begin{equation}\label{defPsi}
\Psi_n(k,x,y)=D_xu_n(k,S_{\Delta t}x+\Delta tS_{\Delta t}\overline{F}(x)).(S_{\Delta t}F(x,y)).
\end{equation}
Then if we define $\overline{\Psi}_{n}(k+1,x)=\int_{H}\Psi_n(k+1,x,y)\mu^x(dy)$ we have
\begin{align*}
|\E D_xu_n(k+1,S_{\Delta t}X_k+\Delta t S_{\Delta t}&\overline{F}(X_k)).(S_{\Delta t}\overline{F}(X_k)-S_{\Delta t}\tilde{F}_k)|\\
&\leq \frac{1}{MN}\sum_{m=n_T}^{n_T+N-1}\sum_{j=1}^{M}|\E \Psi_n(k+1,X_k,Y_{k,m,j})-\E\overline{\Psi}_{n}(k+1,X_k)|.
\end{align*}
By using conditional expectation $\E_k$ with respect to $\mathcal{G}_k$ defined by \eqref{defsigmafield},
$$\E \Psi_n(k+1,X_k,Y_{k,m,j})=\E\E_k \Psi_n(k+1,X_k,Y_{k,m,j})$$
does not depend on $j$, so that in the sequel we fix $j\in\left\{ 1,\ldots,M\right\}$.


The following Lemma gives the regularity results for the auxiliary functions:
\begin{lemme}\label{regPsi_nk}
For any $0<T<+\infty$ and $\Delta t_0>0$, there exists a constant $C$, such that for any $0<\Delta t\leq \Delta t_0$, any $1\leq n\leq n_0=\lfloor \frac{T}{\Delta t}\rfloor$ and any $0\leq k\leq n-1$, the following derivatives exist and are controlled: for any $x,y\in H$, $h,k\in H$,
\begin{gather*}
|D_y\Psi_n(k,x,y).h|\leq C|h|_H,\\
|D_{yy}^{2}\Psi_n(k,x,y).h|\leq C(1+\frac{C}{((n-k)\Delta t)^\eta})|h|_H|k|_H.
\end{gather*}
\end{lemme}
Like in Lemma \ref{regu_n}, the proof relies smoothing effect of the semi-group $(S_{\Delta t}^{n})_{n\in \N}$. The Lemma is proved below in Section \ref{sect_aux_proofs}.

We can then apply Theorem \ref{weakdistance}, valid with the weak dissipativity condition \eqref{hypdiss}, using the conditional expectation $\E_k$: for any $n_T\leq m\leq m_0$ and $k<n-1$
$$|\E \Psi_n(k+1,X_k,Y_{k,m,j})-\E\overline{\Psi}_{n}(k+1,X_k)|\leq Ce^{-cm\tau}+C(1+\frac{1}{((m-1)\tau)^{1/2-\kappa}})\tau^{1/2-\kappa}.$$

Therefore
\begin{align*}
|\E D_xu_n(k+1,S_{\Delta t}X_k+\Delta t S_{\Delta t}&\overline{F}(X_k)).(S_{\Delta t}\overline{F}(X_k)-S_{\Delta t}\tilde{F}_k)|\\
&\leq \frac{1}{MN}\sum_{m=n_T}^{n_T+N-1}\sum_{j=1}^{M}|\E \Psi_n(k+1,X_k,Y_{k,m,j})-\E\overline{\Psi}_{n}(k+1,X_k)|\\
&\leq C\frac{e^{-cn_T\tau}}{N\tau+1}+C\tau^{1/2-\kappa}+C\frac{\tau^{1/2-\kappa}}{((n_T-1)\tau)^{1/2-\kappa}}.
\end{align*}
To conclude, it remains to use \eqref{decompoweak2} and \eqref{thedecompoweak}.

\begin{rem}
When the strict dissipativity Assumption is satisfied, we can indeed obtain a bound without $\frac{1}{((m-1)\tau)^{1/2-\kappa}}$: we can control the distance between the invariant measures of the continuous and discrete time processes, thanks to the second part of Corollary \ref{weakcor}.
\end{rem}

\subsection{Proof of the auxiliary Lemmas \ref{regu_n} and \ref{regPsi_nk}}\label{sect_aux_proofs}

\hspace{100 pt}
\bigskip

\underline{Proof of Lemma \ref{regu_n}}
We use the following expression for $\overline{X}_{n}(x)$:
$$\overline{X}_{n}(x)=S_{\Delta t}^{n}x+\Delta t\sum_{k=0}^{n-1}S_{\Delta t}^{n-k}\overline{F}(\overline{X}_{k}(x)).$$
By definition, for any $0\leq k\leq n$ we have $u_n(k,x)=\Phi(\overline{X}_{n-k}(x))$;
we see that the derivatives in directions $h,h_1,h_2\in H$ are given by
$$D_xu_n(k,x).h=D\Phi(\overline{X}_{n-k}(x)).(\frac{d}{dx}\overline{X}_{n-k}(x).h),$$
and
\begin{align*}
D_{xx}^{2}u_n(k,x).(h_1,h_2)&=D\Phi(\overline{X}_{n-k}(x)).(\frac{d^2}{dx^2}\overline{X}_{n-k}(x).(h_1,h_2))\\
&\hspace{10 pt}+D^2\Phi(\overline{X}_{n-k}(x)).(\frac{d}{dx}\overline{X}_{n-k}(x).h_1,\frac{d}{dx}\overline{X}_{n-k}(x).h_2).
\end{align*}
$\Phi$ is of class $\mathcal{C}^2$ on $H$, with bounded derivatives; therefore we just need to control $\frac{d}{dx}\overline{X}_{n-k}(x).h$ and $\frac{d^2}{dx^2}\overline{X}_{n-k}(x).(h_1,h_2)$.
We use the following estimates of the derivatives of $\overline{F}$, given in Proposition \ref{FbarLip}: for any $x\in H$, $h\in H$, $h_1,h_2 H$,
\begin{gather*}
|D\overline{F}(x).h|\leq C|h|\\
|(-A)^{-\eta}D^2\overline{F}(x).(h_1,h_2)|\leq C|h_1||h_2|.
\end{gather*}

{\em (i)} For any $0\leq n\leq n_0$, we can write
$$\frac{d}{dx}\overline{X}_n(x).h=S_{\Delta t}^{n}h+\Delta t\sum_{k=0}^{n-1}S_{\Delta t}^{n-k}D\overline{F}(\overline{X}_k(x))(\frac{d}{dx}\overline{X}_k(x).h);$$
therefore
$$|\frac{d}{dx}\overline{X}_n(x).h|\leq |h|+C\Delta t\sum_{k=0}^{n-1}|\frac{d}{dx}\overline{X}_k(x).h|,$$
and a discrete Gronwall Lemma then yields
$$|\frac{d}{dx}\overline{X}_n(x).h|\leq |h|e^{Cn\Delta t}\leq |h|e^{CT}.$$

{\em (ii)} For any $0\leq n\leq n_0$, we can write
\begin{align*}
\frac{d^2}{dx^2}\overline{X}_n(x).(h_1,h_2)&=\Delta t\sum_{k=0}^{n-1}S_{\Delta t}^{n-k}(-A)^{\eta}(-A)^{-\eta}D^2\overline{F}(\overline{X}_k(x)).(\frac{d}{dx}\overline{X}_k(x).h_1,\frac{d}{dx}\overline{X}_k(x).h_2)\\
&\hspace{10 pt}+\Delta t\sum_{k=0}^{n-1}S_{\Delta t}^{n-k}D\overline{F}(\overline{X}_k(x)).(\frac{d^2}{dx^2}\overline{X}_k(x).(h_1,h_2)).
\end{align*}
Since $|S_{\Delta t}^{n-k}(-A)^{\eta}|\leq \frac{C}{((n-k)\Delta t)^\eta}$ when $k<n$, and thanks to the previous estimates on $\frac{d}{dx}\overline{X}_k(x).h$ and $D^2\overline{F}$, we get
$$|S_{\Delta t}^{n-k}(-A)^{\eta}(-A)^{-\eta}D^2\overline{F}(\overline{X}_k(x)).(\frac{d}{dx}\overline{X}_k(x).h_1,\frac{d}{dx}\overline{X}_k(x).h_2)|\leq \frac{C}{((n-k)\Delta t)^\eta}|h|_1|h|_2.$$
Therefore
\begin{align*}
|\frac{d^2}{dx^2}\overline{X}_n(x).(h_1,h_2)|&\leq C|h_1||h_2|\\
&\hspace{10 pt}+C\Delta t\sum_{k=0}^{n-1}|\frac{d^2}{dx^2}\overline{X}_k(x).(h_1,h_2)|,
\end{align*}
and a discrete Gronwall Lemma then yields
$$|\frac{d^2}{dx^2}\overline{X}_n(x).(h_1,h_2)|\leq C|h_1||h_2|.$$

{\em (iii)} To prove the last estimate of the Lemma, we write
\begin{align*}
|\frac{d}{dx}\overline{X}_n(x).h|&=|S_{\Delta t}^{n}h+\Delta t\sum_{k=0}^{n-1}S_{\Delta t}^{n-k}D\overline{F}(\overline{X}_k(x))(\frac{d}{dx}\overline{X}_k(x).h)|\\
&\leq \frac{C}{(n\Delta t)^\eta}|h|_{(-A)^{-\eta}}+C\Delta t\sum_{j=0}^{n-1}|\frac{d}{dx}\overline{X}_k(x).h|.
\end{align*}
We obtain
\begin{align*}
(n\Delta t)^\eta|\frac{d}{dx}\overline{X}_n(x).h|&\leq C|h|_{(-A)^{-\eta}}+C\Delta t(n\Delta t)^\eta|h|_{H}\\
&\hspace{10 pt}+C(n\Delta t)^\eta\Delta t\sum_{j=1}^{n-1}\frac{1}{(j\Delta t)^\eta}(j\Delta t)^\eta|\frac{d}{dx}\overline{X}_j(x).h|;
\end{align*}
To conclude, we use Gronwall Lemma, to get
$$(n\Delta t)^\eta|\frac{d}{dx}\overline{X}_n(x).h|\leq C_T(|h|_{(-A)^{-\eta}}+\Delta t(n\Delta t)^\eta).$$\qed

\bigskip
\underline{Proof of Lemma \ref{regPsi_nk}}

{\em (i)} The first derivative with respect to $y$ is easy to control: we have for any $h\in H$
$$D_y\Psi_n(k,x,y).h=D_xu_n(k,S_{\Delta t}x+\Delta tS_{\Delta t}\overline{F}(x)).(S_{\Delta t}D_yF(x,y).h),$$
and we get $|D_y\Psi_n(k,x,y).h|\leq C|h|_H$.

{\em (ii)} When we look at the second-order derivative, we see that
$$D_{yy}^{2}\Psi_n(k,x,y).(h,k)=D_xu_n(k,S_{\Delta t}x+\Delta tS_{\Delta t}\overline{F}(x)).(S_{\Delta t}D_{yy}^{2}F(x,y).(h,k)).$$
Thanks to the last estimate of Lemma \ref{regu_n}, we can control the expression by
$$C\left(\Delta t|S_{\Delta t}D_{yy}^{2}F(x,y).(h,k)|_{H}+\frac{|S_{\Delta t}D_{yy}^{2}F(x,y).(h,k)|_{(-A)^{-\eta}}}{((n-k)\Delta t)^\eta}\right),$$

We then notice that
\begin{align*}
\Delta t|S_{\Delta t}D_{yy}^{2}F(x,y).(h,k)|_{H}&\leq C\Delta t^{1-\eta}|(-A)^{-\eta}D_{yy}^{2}F(x,y).(h,k)|\\
&\leq C|h|_{H}|k|_{H},
\end{align*}
since $\eta<1$ and $\Delta t$ is bounded, and thanks to Assumption \ref{hypF2}; the other part is controlled thanks to Assumption \ref{hypF2}.\qed



\bibliographystyle{plain}

\end{document}